\documentclass[12pt]{article}
\usepackage{amsmath}
\usepackage{amsfonts}
\usepackage{amsthm}

\begin{document}

\def\N{\mathbb{N}}
\def\F{\mathbb{F}}
\def\Z{\mathbb{Z}}
\def\R{\mathbb{R}}
\parindent= 3.em \parskip=5pt
\baselineskip=17pt

\centerline{\bf{QUARTIC POWER SERIES IN $\F_3((T^{-1}))$ }}
\centerline{\bf{WITH BOUNDED PARTIAL QUOTIENTS, II }}
\vskip 0.5 cm
\centerline{by}
\centerline { Domingo Gomez (Universidad de Cantabria)}
\centerline {and Alain Lasjaunias (Universit\'e Bordeaux 1)}
\vskip 0.5 cm
\centerline{\bf{1. Introduction and result}}
\par We are concerned with power series in $1/T$ over a finite field, where $T$ is an indeterminate. If the base field is $\F_q$,
the finite field with $q$ elements, these power series belong to the
field $\F_q((T^{-1}))$, which will be here denoted by $\F(q)$. Thus a
nonzero element of $\F(q)$ is represented by
$$\alpha=\sum_{k\leq k_0}u_kT^k \quad \text{ where } k_0\in \Z, u_k\in
\F_q \quad \text{ and }u_{k_0}\neq 0.$$
An absolute value on this field is defined by $\vert \alpha\vert
=\vert T\vert ^{k_0}$ where $\vert T\vert>1$ is a fixed real
number. We also denote by $\F(q)^{+}$ the subset of power series
$\alpha$ such that $\vert \alpha\vert >1$. We know
that each irrational element $\alpha \in \F(q)^{+}$ can be expanded as an infinite
continued fraction. This is denoted
$$\alpha=[a_1,a_2,\dots,a_n,\dots]\quad \text{ where } a_i\in
\F_q[T]\text{ and }\deg(a_i)>0\text{ for }i\geq 1.$$ 
By truncating this expansion we obtain a rational element, called a convergent
to $\alpha$ and denoted by $x_n/y_n$ for $n\geq 1$. The polynomials
$(x_n)_{n\geq 0}$ and $(y_n)_{n\geq 0}$, called continuants, are both
defined by the same recursion formula : $K_n=a_nK_{n-1}+K_{n-2}$ for $n\geq 2$,
with the initials conditions $x_0=1$ and $x_1=a_1$ or $y_0=0$ and
$y_1=1$. The polynomials $a_i$ are called the partial quotients of the
expansion. For $n\geq 1$, we denote
$\alpha_{n+1}=[a_{n+1},a_{n+2},\dots]$, called the complete quotient, and we have
$$\alpha=[a_1,a_2,\dots,a_n,\alpha_{n+1}]=(x_n\alpha_{n+1}+x_{n-1})/(y_n\alpha_{n+1}+y_{n-1}).$$
The reader may consult [S] for a general account on continued fractions
in power series fields. Throughout this note we are considering the case $q=3$.
\par In a previous article [L1], the second author of this note
investigated the existence of particular power series 
in $\F(3)$, algebraic over $\F_3[T]$, having all partial quotients of
degree 1  in
their continued fraction expansion. A first example of such algebraic power
series appeared in [MR] (p. 401-402). In this article [MR], Mills and
Robbins, developing the pioneer work by Baum and Sweet [BS],
introduced in the general case a particular subset of algebraic power
series. These power series, now called hyperquadratic, are irrational
elements $\alpha$ satisfying an equation $\alpha=f(\alpha^r)$ where
$r$ is a power of the characteristic of the base field and $f$ is a
linear fractional transformation with integer (polynomials in $T$) coefficients.  
\par The theorem which we present in this note is an extended version
of the one presented in [L1]. The proof given here is based on a
method developed in [L2]. Actually, this method could be applied to
obtain other hyperquadratic continued fractions with all partial quotients
of degree 1 in other fields $\F(p)$, for any odd prime $p$. We have the
following :
\newline {\bf{Theorem.}}{\emph{ Let $m\in\mathbb{N}^*$,  
    $\boldsymbol{\eta}=(\eta_1,\eta_2,\dots,\eta_m)\in (\F_3^*)^m$ where
    $\eta_m=(-1)^{m-1}$ and  {\bf k}$=(k_1,k_2,\dots,k_m)\in \N^m$ where
    $k_1\ge 2$ and $k_{i+1}-k_i\ge 2 $ for $i=1,\ldots,m-1.$
    We define the following integers,
    $$
      t_{i,n}=k_m(3^n-1)/2+k_i 3^n \quad \text{ for }\quad 1\le i\le m
      \quad \text{ and }\quad  n\ge 0.   $$
    We define two sequences $ (\lambda_t)_{t\ge 1} $ and
    $(\mu_t)_{t\ge 1}$ in $\F_3$. For $n\geq 0$, we have  
    \begin{equation*}
      \lambda_t=
      \begin{cases}
        1 &\text{if } \quad 1\le t\le t_{1,0},\\
        (-1)^{mn+i}& \text{if } \quad t_{i,n}< t\le t_{i+1,n} \quad \text{
          for }1\le i< m,\\
        (-1)^{m(n+1)}& \text{if } \quad t_{m,n}< t\le t_{1,n+1}.\\
      \end{cases}
    \end{equation*}
    Also $\mu_1=1$ and for $n\geq 0$, $1\leq i\leq m$ and $t>1$
    \begin{equation*}
      \mu_t=
      \begin{cases}
        (-1)^{n(m+1)}\eta_i& \text{if }  t= t_{i,n}\text{ or }t= t_{i,n}+1,\\
        0 &\text{otherwise.}
      \end{cases}
    \end{equation*}    
    Let $\omega(m,\boldsymbol{\eta},\text{\bf k})\in\F(3)$ be defined by
    the infinite continued fraction expansion
    $$\omega(m, \boldsymbol{\eta},\text{\bf k})=
    [a_1,a_2,\ldots,a_n,\ldots]\quad \text{ where }\quad  a_n=\lambda_n T+\mu_n.$$
    We set $l=1+k_m$ and we consider the two usual
    sequences $(x_n)_{n\geq 0}$ and $(y_n)_{n\geq 0}$ as being the
numerators and denominators of the convergents to
$\omega(m,\boldsymbol{\eta},\text{\bf k})$. 
\newline Then $\omega(m,\boldsymbol{\eta},\text{\bf k})$ is the unique root in $\F(3)^+$ of
the quartic equation
$$ X=\frac{x_lX^3+(-1)^{m-1}x_{l-3}}{y_lX^3+(-1)^{m-1}y_{l-3}}.$$}}
\vskip 0.5 cm
\noindent{\bf{Remark.}}{\emph{ The case $m=1$ and thus $\boldsymbol{\eta}=(1)$, {\bf
      k}\,$=(k_1)$, of this theorem is proved in [L1]. The case $m=2$, $\boldsymbol{\eta}=(-1,-1)$ and $\text{\bf k}=(3,6)$
corresponds to the example which was introduced by W. Mills and D. Robbins,
in [MR].}} 
\vskip 0.5 cm
\centerline{\bf{2. Proof of the theorem}}
\par The proof will be divided into three steps.
\newline $\bullet$ {\emph{First step of the proof:}} According to [L2,
Theorem 1, p. 332], there exists a unique element $\beta \in \F(3)$ defined by
$\beta=[a_1,\dots,a_l,\beta_{l+1}]$  and satisfying   
$$\beta^3=(-1)^m(T^2+1)\beta_{l+1}+T+1\quad \text{ and
}\quad a_i=\lambda_i T+ \mu_i,\text{ for }1\le i\le l,$$
where $\lambda_i,\ \mu_i$ are the elements defined in the theorem.
We know that this element is hyperquadratic and that it is the unique root 
in $\F(3)^+$ of the algebraic
equation $X=(x_lX^3+B)/(y_lX^3+D)$ where
$$B=(-1)^m(T^2+1)x_{l-1}-(T+1)x_l \quad \text{ and
}\quad D=(-1)^m(T^2+1)y_{l-1}-(T+1)y_l.$$
We need to transform $B$ and $D$. Using the recursive
formulas for the continuants, we can write
$$K_{l-3}=(a_la_{l-1}+1)K_{l-1}-a_{l-1}K_l. \eqno{(1)}$$
The $l$ first partial quotients of
$\beta$ are given, from the hypothesis of the theorem, and we have
$$a_{l-1}=(-1)^{m-1}(T+1)\quad \text{ and
}\quad a_{l}=(-1)^{m-1}(-T+1).\eqno{(2)}$$
Combining (1), applied to both sequences $x$ and $y$, and (2), we get 
$$B=(-1)^{m-1}x_{l-3}\quad \text{ and
}\quad D=(-1)^{m-1}y_{l-3}.$$
Hence we see that $\beta$ is the unique root in $\F(3)^+$ of
the quartic equation stated in the theorem. 
\vskip 0.5 cm
\noindent $\bullet$ {\emph{Second step of the proof:}} In this section
$l\geq 1$ is a given integer. We consider all the infinite continued fractions $\alpha \in \F(3)$ defined by
$\alpha=[a_1,\dots,a_l,\alpha_{l+1}]$ where
$$(I) \quad :\quad a_i=\lambda_iT+\mu_i \quad \text{ with}\quad (\lambda_i,\mu_i)\in
\F_3^*\times \F_3 ,\quad \text{ for}\quad 1\leq i \leq l$$
and  
$$(II) \quad :\quad \alpha^3=\epsilon_1(T^2+1)\alpha_{l+1}+\epsilon_2T+\nu_0 \quad \text{ with}\quad (\epsilon_1,\epsilon_2,\nu_0)\in
\F_3^*\times \F_3^*\times \F_3.$$
See [L2, Theorem 1, p. 332], for the existence and unicity of $\alpha
\in \F(3)$ defined by
the above relations. Our aim is to show that these continued fraction
expansions can be explicitly described, under particular conditions on the
parameters $(\lambda_i,\mu_i)_{1\leq i\leq l}$ and
$(\epsilon_1,\epsilon_2,\nu_0)$. Following the same method as in [L2],
we first prove :  
\newline {\bf{Lemma 1.}}{\emph{ Let $(\lambda, \epsilon_1,\epsilon_2)
    \in (\F_3^*)^3$ and $\nu \in \F_3$. We set $U=\lambda T^3-\epsilon_2 T+\nu$,
     and $V=\epsilon_1(T^2+1)$. We set $\delta=\lambda
    +\epsilon_2$ and we assume that $\delta\neq 0$. We define $ \epsilon^*=1$
    if $\nu =0$ and $\epsilon^*=-1$ if $\nu \neq 0$. Then the continued fraction
     expansion for $U/V$ is given by  
$$ U/V=[\epsilon_1\lambda T,-\epsilon_1(\delta T+\nu)
,-\epsilon_1(\epsilon^*\delta T+\nu)].$$
Moreover, setting $U/V=[u_1,u_2,u_3]$, then for $X\in \F(3)$
we have  $$[U/V,X]=[u_1,u_2,u_3,\frac{X}{(T^2+1)^2}+\frac{\epsilon^*\epsilon_1(\delta
  T+\nu)}{T^2+1}].$$}}
Proof : Since $\epsilon_1^2=1$ and $\delta^2=1$, we can write
$$U=\epsilon_1\lambda TV-\delta T+\nu \quad \text{ and}\quad
V=\epsilon_1(\delta T+\nu)(\delta
T-\nu)+\epsilon_1(1+\nu^2).\eqno{(3)}$$
Clearly (3) implies the following continued fraction expansion
$$ U/V=[\epsilon_1\lambda T,-\epsilon_1(\delta
T+\nu),\epsilon_1(1+\nu^2)(-\delta T+\nu)].\eqno{(4)}$$
Finally, observing that $\epsilon_1(1+\nu^2)=\epsilon^*\epsilon_1$ and
$\epsilon^*\epsilon_1\nu=-\epsilon_1\nu$, we see that (4) is the
expansion stated in the lemma. The last formula is obtained from [L2] (Lemma
3.1 p. 336). According to this lemma, we have
$$[U/V,X]=[u_1,u_2,u_3,X']\quad \text{ where}\quad
X'=X(u_2u_3+1)^{-2}-u_2(u_2u_3+1)^{-1}.$$
We check that $u_2u_3=T^2$ if $\nu =0$ and  $u_2u_3=\nu^2-T^2$ if $\nu
\neq 0$, therefore we have $u_2u_3+1=\epsilon^*(T^2+1)$ and this
implies the desired equality.
\newline We shall prove now a second lemma. In the sequel we define $f(n)$ as $3n+l-2$ for $n\geq 1$. We have the
following :
\newline {\bf{Lemma 2.}}{\emph{ Let $\alpha=[a_1,\dots,a_n,\dots]$ be
    an irrational element of $\F(3)$. We assume that for an index $n\geq 1$ we
have 
$a_n=\lambda_nT+\mu_n$ with $(\lambda_n,\mu_n)\in
\F_3^*\times \F_3$ and 
$$\alpha_n^3=\epsilon_1(T^2+1)\alpha_{f(n)}+\epsilon_{2,n}T+\nu_{n-1}\quad
\text{ where }\quad
(\epsilon_1,\epsilon_{2,n},\nu_{n-1})\in
(\F_3^*)^2\times \F_3.$$
 We set $\nu_{n}=\mu_n-\nu_{n-1}$ and
$\epsilon_n^*=1$ if $\nu_{n}=0$ or $\epsilon_n^*=-1$ if
$\nu_{n}\neq 0$. We set $\delta_n=\lambda_n+\epsilon_{2,n}$, and
$\epsilon_{2,n+1}=-\epsilon_n^*\delta_n$. We assume that
$\delta_n\neq 0$. Then we have :
$$(a_{f(n)},a_{f(n)+1},a_{f(n)+2})=(\epsilon_1\lambda_nT,-\epsilon_1(\delta_nT+\nu_{n}),-\epsilon_1(\epsilon_n^*\delta_nT+\nu_{n}))$$
and
$$\alpha_{n+1}^3=\epsilon_1(T^2+1)\alpha_{f(n+1)}+\epsilon_{2,n+1}T+\nu_{n}.$$
}}
Proof: We can write $\alpha_n^3=[a_n^3,\alpha_{n+1}^3]=[\lambda_n
T^3+\mu_n,\alpha_{n+1}^3]$. Consequently 
$$\alpha_n^3=\epsilon_1(T^2+1)\alpha_{f(n)}+\epsilon_{2,n}T+\nu_{n-1}$$ 
 is equivalent to 
$$[(\lambda_nT^3+\mu_n-\epsilon_{2,n}T-\nu_{n-1})/(\epsilon_1(T^2+1)),\epsilon_1(T^2+1)\alpha_{n+1}^3]=\alpha_{f(n)}.\eqno{(5)}$$
Now we apply Lemma 1 with $U=\lambda_nT^3-\epsilon_{2,n}T+\nu_{n}$ and
$X=\epsilon_1(T^2+1)\alpha_{n+1}^3$. Consequently (5) can be written
as
$$[\epsilon_1\lambda_n T,-\epsilon_1(\delta_n T+\nu_n)
,-\epsilon_1(\epsilon_n^*\delta_n
T+\nu_n),X']=\alpha_{f(n)}\eqno{(6)}$$
where
$$X'=(\epsilon_1\alpha_{n+1}^3+\epsilon_1\epsilon_n^*(\delta_n
T+\nu_n))/(T^2+1).\eqno{(7)}$$
Moreover we have $\vert \alpha_{n+1}^3\vert \geq \vert T^3\vert$ and
consequently $\vert X'\vert >1$. Thus (6) implies that the three
partial quotients $a_{f(n)}$, $a_{f(n)+1}$ and $a_{f(n)+2}$ are as
stated in this lemma and also that we have
$X'=\alpha_{f(n+1)}$. Combining this last equality with (7), and
observing that $-\epsilon_n^*\nu_n=\nu_n$, we obtain
the desired statement. 
\par Applying Lemma 2, we see that for a continued fraction defined by
(I) and (II), the partial quotients, from the rank $l+1$ onward, can
be given explicitly three by three, as long as the quantity $\delta_n$ is not zero. This is
taken up in the following proposition :
\newline {\bf{Proposition 3.}}{\emph{ Let $\alpha \in \F(3)$ be an
    infinite continued fraction expansion defined by (I) and (II).
Then there exists $N\in \N^*\cup \lbrace \infty \rbrace$ such that
$$a_n=\lambda_nT+\mu_n \quad \text{ where } \quad (\lambda_n,\mu_n)\in
\F_3^*\times \F_3 \quad   \quad \text{ for }\quad 1\leq n < f(N).$$ 
For $1\leq n <f(N)$, we define $\nu_n=\sum_{1\leq i\leq
  n}(-1)^{n-i}\mu_i+(-1)^n\nu_0$. Then we have 
$$\mu_{f(n)}=0 \quad
\text{ and  }\quad \mu_{f(n)+1}=\mu_{f(n)+2}=-\epsilon_1\nu_n \quad
\text{ for }\quad 1\leq n < N.$$
For $1\leq n <N$, we define $\epsilon_n^*=1$ if $\nu_{n}=0$ or $\epsilon_n^*=-1$ if
$\nu_{n}\neq 0$. We define recursively the sequence $(\delta_n)_{1\leq n\leq N}$ by  
$$\delta_1=\lambda_1+\epsilon_2 \quad
\text{ and  }\quad \delta_n=\lambda_n-\epsilon_{n-1}^*\delta_{n-1}\quad
\text{ for  }\quad 2\leq n\leq N.$$
Then, for $1\leq n < N$, we have 
$$\lambda_{f(n)}=\epsilon_1\lambda_n, \quad
\lambda_{f(n)+1}=-\epsilon_1\delta_n \quad
\text{ and  }\quad
\lambda_{f(n)+2}=-\epsilon_1\epsilon_n^*\delta_n.$$}}
Proof: Starting from (II), since $f(1)=l+1$, setting
$\epsilon_2=\epsilon_{2,1}$ and observing that all the partial
quotients are of degree 1, we can apply repeatedly
Lemma 2 as long as we have $\delta_n\neq 0$. If $\delta_n$ happens to
vanish, the process is stopped and we denote by $N$ the first index such
that $\delta_N=0$, otherwise $N$ is $\infty$. From the formulas
$\delta_n=\lambda_n+\epsilon_{2,n}$ and $\epsilon_{2,n+1}=-\epsilon_n^*\delta_n$
for $n\geq 1$, we obtain the recursive formulas for the sequence
$\delta$. Moreover, the formula $\nu_n=\mu_n-\nu_{n-1}$, implies clearly the
equality for $\nu_n$. Finally the formulas concerning $\mu$ and $\lambda$ are
directly derived from the three partial quotients $a_{f(n)}$,
$a_{f(n)+1}$ and $a_{f(n)+2}$ given in Lemma 2.

\vskip 0.5 cm 
\noindent $\bullet$ {\emph{Last step of the proof:}}
We start from the element $\beta\in\F(3)$, introduced in 
the first step of the proof, defined  
by its $l$ first partial quotients, where $l=k_m+1$, and by $(II)$ with
$(\epsilon_1,\epsilon_2,\nu_0)=((-1)^m,1,1)$. 
According to the first step of the proof, we 
need to show that $\beta=\omega(m,\boldsymbol{\eta},\text{\bf k})$. To
do so,  we apply Proposition 3 to $\beta$, and we show that $N=\infty$ and that
the resulting sequences $(\lambda_n)_{n\ge 1}$
and $(\mu_n)_{n\ge 1}$ are the one which are described in the theorem. 
\newline From the $l$-tuple $(\mu_1,\dots,\mu_l)$ and $\nu_o=1$, we
obtain 
$$\nu_{t}=\eta_i \quad
\text{ if }\quad t=t_{i,0} \quad
\text{ and }\quad \nu_{t}=0 \quad
\text{ otherwise, }\quad
\text{ for }\quad 1\leq t\leq
l.\eqno{(8)}$$
Since $\mu_{f(n)+1}=\mu_{f(n)+2}$, we have
$\nu_{f(n)+2}=\nu_{f(n)}$. Since $\mu_{f(n)}=0$, we also have
$\nu_{f(n)}=-\nu_{f(n)-1}=-\nu_{f(n-1)+2}$. This implies $\nu_{f(n)+2}=(-1)^{n-1}\nu_{f(1)+2}$. Since $\nu_{f(1)+2}=\nu_{f(1)}=-\nu_{f(1)-1} =-\nu_{l}=0$, we obtain
  $$  \nu_{f(n)}=\nu_{f(n)+2}=0\quad \text{for }\quad 1\le n < N. \eqno{(9)}
 $$
Moreover, from $\nu_{f(n)+1}=\mu_{f(n)+1}-\nu_{f(n)}$ and (9), we also
get
  $$  \nu_{f(n)+1}=-\epsilon_1\nu_n \quad \text{for }\quad 1\le n < N. \eqno{(10)}
 $$
Now, it is easy to check that we have $f(t_{i,n})+1=t_{i,n+1}$.
Since $\epsilon_1=(-1)^m$, (10) implies
$\nu_{t_{i,n}}=(-1)^{m+1}\nu_{t_{i,n-1}}$ if $t_{i,n}<f(N)$. By
induction from (8), with (9) and (10), we obtain 
$$\nu_{t_{i,n}}=(-1)^{(m+1)n}\eta_i \quad
\text{ and }\quad \nu_{t}=0 \text{ if } t\neq t_{i,n}, \text{ for }\quad 1\leq t <f(N). \eqno{(11)}$$
Since we have $\mu_{n}=\nu_{n}+\nu_{n-1}$, from (11) and $\nu_0=1$, we see that
$\mu_{n}$ satisfies the formulas given in the theorem, for $1\leq n <f(N)$. 
Moreover, (11) implies clearly the following :
$$\epsilon^*_t=
      \begin{cases}
        -1& \text{if }  t= t_{i,n},\\
        1 &\text{otherwise}
      \end{cases} \quad \text{ for }\quad 1\leq t <f(N). \eqno{(12)} $$   
Now we turn to the definition of the sequence $(\lambda_n)_{n\geq 1}$
given in the theorem, corresponding to the element $\omega$. With our notations and according to (12), we observe that this definition can be translated into the following formulas 
$$\lambda_1=1 \quad \text{ and }\quad
\lambda_{n}=\epsilon_{n-1}^*\lambda_{n-1}\quad \text{ for } \quad 2\leq n
< f(N). \eqno{(13)}$$
Consequently, to complete the proof, we need to establish that $N=\infty$ and that
(13) holds. The recurrence relation binding the sequences
$\delta$ and $\lambda$ , introduced in Proposition
3, can be written as
$$\delta_n+\lambda_n=-\epsilon_{n-1}^*(\delta_{n-1}+\lambda_{n-1})+\epsilon_{n-1}^*\lambda_{n-1}-\lambda_n\quad
\text{ for } \quad 2\leq n \leq N. \eqno{(14)}$$
Comparing (13) and (14), we see that $\delta_n+\lambda_n=0$, for $n\geq 1$, will
imply that $\delta_n$ never vanishes, i.e. $N=\infty$, and that the sequence
$(\lambda_n)_{n\geq 1}$ is the one which is described in the theorem. So we only
need to prove that $\delta=-\lambda$. Since $\beta$ and $\omega$ have
the same first partial quotients, (13) holds for $2\leq n\leq l$. Since $\delta_1=\lambda_1+\epsilon_2=-1=-\lambda_1$, combining (13)
and (14), we obtain $\delta_n=-\lambda_n$ for $1\leq n\leq l$. We
also have
$\lambda_{l+1}=\lambda_{f(1)}=\epsilon_1\lambda_1=(-1)^m=\lambda_l$,
and therefore we get
$\delta_{l+1}=\lambda_{l+1}-\epsilon_l^*\delta_l=\lambda_{l+1}+\lambda_{l}=-\lambda_{l+1}$. By
induction, we
shall now prove that $\delta_t=-\lambda_t$ for
$t=f(n)+1,f(n)+2$ anf $f(n+1)$ with $n\geq 1$. From (9) and (10), we
have
$ \epsilon_{f(n)}^*=\epsilon_{f(n)+2}^*=1$ and
$\epsilon_{f(n)+1}^*=\epsilon_n^*$. Thus we get
$$
\delta_{f(n)+1}=\lambda_{f(n)+1}-\epsilon_{f(n)}^*\delta_{f(n)}=\lambda_{f(n)+1}+ \lambda_{f(n)}=-\epsilon_1\delta_n+\epsilon_1\lambda_n=-\lambda_{f(n)+1}.
$$
$$
\delta_{f(n)+2}=\lambda_{f(n)+2}-\epsilon_{f(n)+1}^*\delta_{f(n)+1}=\lambda_{f(n)+2}+\epsilon_n^*\lambda_{f(n)+1}=-\lambda_{f(n)+2}.$$
 $$ \delta_{f(n+1)}=\lambda_{f(n+1)}-\epsilon_{f(n)+2}^*\delta_{f(n)+2}=\epsilon_1\lambda_{n+1}+\lambda_{f(n)+2}=\epsilon_1(\lambda_{n+1}-\epsilon_n^*\delta_n)$$
 $$ = \epsilon_1\delta_{n+1}=-\epsilon_1\lambda_{n+1}=-\lambda_{f(n+1)}.$$
So the proof of the theorem is complete.
\par Finally, we make a remark concerning the apparent peculiarity of
the theorem which is proved in this note. We make the following
conjecture: Let $\alpha \in \F(3)$ be a hyperquadratic element, which
is not quadratic; then
$\alpha$ has all its partial quotients of degree 1, except for the
first ones, if and only if there exist a linear fractional
transformation
 $f$, with coefficients in $\F_3[T]$ and determinant in
$\F_3^*$ , a triple $(m,\boldsymbol{\eta},\text{\bf k})$ and a pair $(\lambda,\mu)\in
\F_3^*\times \F_3$ such that
$\alpha(T)=f(\omega(m,\boldsymbol{\eta},\text{\bf k})(\lambda T+\mu))$.

\vskip 1 cm

\begin{tabular}[b]{|l|r|}
    \multicolumn{2}{}{} \\
   Domingo Gomez & Alain Lasjaunias \\
   Universidad de Cantabria & Universit\'e Bordeaux 1 \\ 
   Departamento de Matem\'aticas &  C.N.R.S.-UMR 5251 \\
   39005 Santander, Spain & 33405 Talence, France \\
   Domingo.Gomez@unican.es & Alain.Lasjaunias@math.u-bordeaux1.fr \\
\end{tabular}

\end{document}